\documentclass[11pt]{article}
\usepackage{amssymb,latexsym,amsmath,amsbsy,amsthm,amsxtra,amsgen}
\newfont{\msbm}{msbm10 at 11pt}

\newcommand {\Z} {\mbox{\msbm Z}}

\def\mn{\bigskip\noindent}

\newtheorem{Theo}{Theorem}
\newtheorem{Lemma}[Theo]{Lemma}

\newtheorem{Rmk}[Theo]{Remark}

\begin{document}
\title{A contact process with mutations on a tree}
\author{by Thomas M. Liggett\thanks{Partially supported by NSF grant DMS-0301795}, \:Rinaldo B. Schinazi, \:Jason Schweinsberg\thanks{Partially supported by NSF grant DMS-0504882}\\\\
University of California at Los Angeles,\\
 University of Colorado at Colorado Springs \\
and
 University of California at San Diego}
\maketitle

\begin{abstract}
Consider the following stochastic model for immune response.
Each pathogen gives birth to a new pathogen at rate $\lambda$.  When a new pathogen is born, it has the same type as its parent with probability $1 - r$.  With probability $r$, a mutation occurs, and the new pathogen has a different type from all previously observed pathogens. When a new type appears in the population, it survives for an exponential amount of time with mean $1$, independently of all the other
types.  All pathogens of that type are killed simultaneously. Schinazi and Schweinsberg (2006) have shown that this model on $\Z^d$ behaves rather differently from its non-spatial version. In this paper, we show that this model on a homogeneous tree captures features from both the non-spatial version and the $\Z^d$ version.  We also obtain comparison results between this model and the basic contact process on general graphs.

\end{abstract}

\footnote{{\it Key words and phrases}:  mutation, immune system, branching process, spatial stochastic model, contact process}
\footnote{{\it 2000 Mathematics Subject Classification}: 60K35}

\section{Introduction}

Schinazi and Schweinsberg (2006) have recently introduced a non-spatial and a spatial version of the following stochastic model for immune response. For the non-spatial version (Model 2 in their paper), each pathogen gives birth to a new pathogen at rate $\lambda$.  When a new pathogen is born, it has the same type as its parent with probability $1 - r$.  With probability $r$, a mutation occurs, and the new pathogen has a different type from all previously observed pathogens. When a new type appears in the population, it survives for an exponential amount of time with mean $1$, independently of all the other
types.  All pathogens of that type are killed simultaneously. 

For the spatial version (Model S2 in Schinazi and Schweinsberg (2006)), let $S$ be the square lattice 
$\Z^d$ or the homogeneous tree $T_d$ in which every vertex has $d+1$ neighbors.  Let $x$ be a site in $S$ occupied by a pathogen and $y$ be one of its nearest neighbors. There are $2d$ such neighbors in $\Z^d$ and $d+1$ in $T_d$. After a random exponential time with rate $\lambda$, the pathogen on $x$ gives birth on $y$, provided $y$ is empty (if $y$ is occupied nothing happens). With probability $1-r$ the new pathogen on $y$ is of the same type as the parent pathogen on $x$. With probability $r$ the new pathogen is of a different type. We assume that every new type that appears is different from all types that have ever appeared.
As in the non-spatial model, we assume that each type survives for an exponential amount of time with mean $1$, independently of all other types, and that all pathogens of a type are killed simultaneously.

Biological motivation has been provided in Schinazi and Schweinsberg (2006). There is also some obvious motivation for such a model for the spread of a virus in computer networks. In particular, when a virus is discovered in a computer network it is usually destroyed at once. It is also interesting to test the influence of particular network topologies on the spread of a virus;
see for instance Pastor-Sattoras and Vespignani (2001), Berger et. al. (2005), Ganesh et. al. (2005), and Draief et al. (2006).

Schinazi and Schweinsberg (2006) have shown that the pathogens survive in the non-spatial model
with positive probability if and only if $\lambda>1$ and $r>0$; see their Theorem 1.2. In contrast, for the spatial model on $\Z^d$ they have shown that if $\lambda$ is large enough then the pathogens survive
with positive probability for large $r$ and die out with probability 1 for small $r$; see their Theorem 4. Given that a ball of radius $R$ has of the order of $d^R$ sites in the homogeneous tree $T_d$ one
might conjecture that the behavior of the spatial model on $T_d$ is similar to that of the non-spatial model (there are so many sites on the tree that space may not be a limitation). As the reader will see in 
the next section, this is not so. In fact, the model on the tree captures both features from the model on $\Z^d$ and from the non-spatial model. These results will be stated in Section 2 and proved in Sections 4 and 5.

The basic contact process may be thought of as being a particular case of this model with $r=1$. See Liggett (1999) for background and results on the contact process. For $r=1$, all pathogens are of 
different types and therefore only one pathogen dies at a time. 
Hence, one might conjecture that the model with $r=1$ has a better chance of surviving than any model 
with $r<1$. This turns out to be true but requires nontrivial arguments. This will be done in Section 3.
This comparison result applies to any graph.

\section{Phase transitions}

We say that the pathogens survive on a graph $S$ if there is a positive probability that at all times there is at least one pathogen somewhere in $S$. If the pathogens do not survive they are said to die out.

\begin{Theo} 
Consider the contact process with mutations on the homogeneous tree $T_d$ with $d \geq 2$, started with a single pathogen of type $1$.
\begin{enumerate}

\item If $\lambda> \frac{1}{d-1}$ then the pathogens survive for all $r>0$. 

\item If $\lambda\leq \frac{1}{d-1+2r}$ then the pathogens die out.  In particular, if $\lambda \leq \frac{1}{d+1}$, the pathogens die out for all $r \geq 0$.

\item If $\frac{1-d+\sqrt{(d-1)(7+9d)}}{2(d^2-1)}<\lambda <\frac 1{d-1}$ the pathogens die out for $r$ close to 0 and survive for $r$ close to 1.
\end{enumerate}
\label{phasetransition}
\end{Theo}

\noindent These results will be proved in Sections 4 and 5.

Note that for all $d \geq 2$, the left-hand side of the inequality in Theorem 1.3 is strictly smaller than the right-hand side.  Therefore, Theorem 1.3 shows that for intermediate values of $\lambda$ there is a phase transition in $r$. Moreover, the difference between the right and left sides of Theorem 1.3 is asymptotic to $\frac 2{3d^2}$ as $d\uparrow\infty$. Combining Theorem 1.1 and 1.2 , we see that the
set of $\lambda$'s for which there is a phase transition in $r$ is of size asymptotically
at most $\frac 2{d^2}$. So, up to a factor of three, we have found the correct size of this set.

We now turn to another type of phase transition. Let $A_t$ be the set of sites that are occupied by any pathogen. The process is started with a single pathogen at the site $x$.
Recall that pathogens are said to survive if
\begin{equation}
P(A_t\not=\emptyset,\forall t>0)>0.
\label{survivecond}
\end{equation}
The pathogens are said to survive weakly if (\ref{survivecond}) holds and
\begin{equation}
P(x \notin A_t \mbox{ for sufficiently large }t) = 1.
\label{weakcond}
\end{equation}
That is, with positive probability, there are pathogens somewhere in $S$ for all times, but almost surely a given site gets infected only finitely many times.  If $P(x \in A_t \mbox{ for arbitrarily large }t) > 0$, then the process is said to survive strongly.

The interest in interacting particles on trees was sparked by Pemantle (1992). He proved that, in contrast to what happens on $\Z^d$, there exists an intermediate range of $\lambda$'s for which the
basic contact process on a homogeneous tree survives weakly; see Liggett (1999) for a complete description of the results and people involved. We next give a partial result in the same direction for the contact process with mutations.

\begin{Theo} Consider the contact process with mutations on the homogeneous tree $T_d$.
If $\frac{1}{d-1}<\lambda< \frac{1}{2\sqrt{d}}$ then the pathogens survive weakly for all $r>0$. 
\label{localdeath}
\end{Theo}
This proves that for $d \geq 6$, there is a range of $\lambda$ for which the contact process with mutations survives weakly. We believe this result to hold for all $d\geq 2$. 
Theorem 2 will be proved at the end of Section 3.

\section{\bf Comparison between the contact process with mutations and the basic contact process}

Let $S$ be a connected undirected graph.
Let $A_t$ be the set of sites occupied by a pathogen at time $t$ for the contact process with mutations on $S$. A typical state is a finite subset $A$ of $S$, together with a partition
$\{A^i\}$ of $A$; $A^i$ is the set of sites of type $i$. Denote the semigroup and
generator of the process by $S_r(t)$ and
$L_r$. Note that the mutation model with
$r=1$ can be viewed as a version of the basic contact process, but only if all initial
types are distinct; otherwise they are not the same. To make the comparison with the
basic contact process, it is convenient to make it have the same birth mechanism as
the multitype version. We will carry along the information about types, but ignore that
information in the death mechanism -- each individual dies separately at rate 1. We
will use primes when referring to this process. Thus the generators are given by
$$L_rg(A)=\sum_i[g(A\backslash A^i)-g(A)]+\lambda\sum_i\sum_{x\in A^i,y\notin A,x\sim y}
[rg(U^yA)+(1-r)g(V_i^yA)-g(A)]$$
and
$$L_r'g(A)=\sum_{x\in A}[g(A\backslash \{x\})-g(A)]+\lambda\sum_i\sum_{x\in A^i,y\notin A,x\sim y}
[rg(U^yA)+(1-r)g(V_i^yA)-g(A)]$$
where $U^yA=V_i^yA=A\cup\{y\}$ as sets, but with different partitions: In $U^yA$, $\{y\}$
is added as a new set in the partition, while in $V_i^yA$, $y$ is added to $A^i$.
Here $x\sim y$ means thar $x$ and $y$ are neighbors in $S$.

We begin with a result about the basic contact process, which is an extension of Theorem 6.2
of Harris (1974). His result is the following statement in the special case $S=Z^d$ and
$C=Z^d$. The proof below is an adaptation of the proof of Proposition 5.9 of Chapter III
of Liggett (1985).
(As a historical note, we point out that Harris proved his result using
coupling rather than duality. He didn't discover contact process duality until two years later.)
Here $P^A$ refers to probabilities when $A$ is the initial state of the process at time zero.

\begin{Lemma} For all $A,B, C\subset S$ and every $t>0$,
$$P^{A\cap B}(A_t'\cap C\neq\emptyset)+P^{A\cup B}(A_t'\cap C\neq\emptyset)\leq
P^{A}(A_t'\cap C\neq\emptyset)+P^{B}(A_t'\cap C\neq\emptyset).$$
\end{Lemma}

\begin{proof}
By checking cases, one sees that
$$1_{\{C\cap(A\cap B)\neq\emptyset\}}+1_{\{C\cap(A\cup B)\neq\emptyset\}}\leq
1_{\{C\cap A\neq\emptyset\}}+1_{\{C\cap B\neq\emptyset\}}$$
for all $A,B,C$. Replacing $C$ by $A_t'$ and taking expectations gives
$$P^C(A_t'\cap(A\cap B)\neq\emptyset)+P^C(A_t'\cap(A\cup B)\neq\emptyset)\leq
P^C(A_t'\cap A\neq\emptyset)+P^C(A_t'\cap B\neq\emptyset).$$
Now use duality (Theorem 1.7 of Chapter VI of Liggett (1985)).
\end{proof}

Here is the basic comparison result:
\begin {Theo} For every $A, C\subset S$, $0\leq r\leq 1$, and $t\geq 0$,
$$P^A(A_t\cap C\neq\emptyset)\leq P^A(A_t'\cap C\neq\emptyset).$$
\end{Theo}

\begin {proof} Let $f(A)=1_{\{A\cap C\neq\emptyset\}}$. The integration by parts formula gives
\begin{align}
P^A(A_t'\cap C\neq\emptyset)-P^A(A_t\cap C\neq\emptyset)&=S_r'(t)f(A)-S_r(t)f(A) \nonumber \\
&=\int_0^tS_r(s)\big(L_r'-L_r)S_r'(t-s)f(A)ds. \nonumber
\end{align}

The easiest way to check this identity is to integrate the derivative with respect to $s$
of $S_r(s)S_r'(t-s)f$ from 0 to $t$, recalling that while the two semigroups do not
commute, each semigroup commutes with its own generator.
By the lemma, $g(A)=S_r'(t-s)f(A)$ (note that $g(A)$ depends on
$A$ only through the set of occupied sites -- not on its partition by type) satisfies
\begin{equation}
g(A\cap B)+g(A\cup B)\leq g(A)+ g(B)
\label{gAB}
\end{equation}
for all $A,B$. Thus it suffices to show that $(L_r'-L_r)g\geq 0$ for all $g$
that depend on $A$ only through the set of occupied sites and satisfy
(\ref{gAB}). Of course, $L_rg(A)$ does depend on the types.
So, write $A=\cup_i A^i$, where $\{A^i\}$ is the partition of $A$ according to
type. Then since the terms corresponding to births agree for the two generators,

\begin{eqnarray*}
(L_r'-L_r)g(A)&=&\sum_{x\in A}[g(A\backslash\{x\})-g(A)]-
\sum_i[g(A\backslash A^i)-g(A)]\\
&=&\sum_i\bigg[\sum_{x\in A^i}[g(A\backslash\{x\})-g(A)]
-g(A\backslash A^i)+g(A)\bigg].
\end{eqnarray*}
We will check that this is nonnegative by checking that each summand (in the sum on $i$)
is nonnegative. So, we need to check that if $B\subset A$, then
$$\sum_{x\in B}g(A\backslash\{x\})\geq g(A\backslash B)+(|B|-1)g(A).$$
To deduce this from (\ref{gAB}), write $B=\{u_1,...,u_k\}$. Then (\ref{gAB}) implies that
\begin{equation}
g(A\backslash\{u_j\})\geq g(A\backslash\{u_1,...,u_j\})-g(A\backslash\{u_1,...,u_{j-1}\})
+g(A)
\label{setg}
\end{equation}
for each $j$. Now sum (\ref{setg}) for $1\leq j\leq k$ to complete the proof.
\end {proof}

\begin{Rmk}
{\em By taking $C = S$ in Theorem 4, we get that for any graph $S$ and any $\lambda > 0$, if the basic contact process dies out, then the contact process with mutations dies out for every $r\in [0,1]$.}
\end{Rmk}

\begin{proof}[Proof of Theorem 2.] 
By Theorem 1.1, the contact process with mutations survives for $\lambda > \frac{1}{d-1}$.  Therefore, it remains only to verify (\ref{weakcond}) when $\lambda < \frac{1}{2 \sqrt{d}}$.  Recall that $x$ is the 
location of the pathogen at time zero.  Following Section 1.4 of Liggett (1999), define a function $\ell: T_d \rightarrow \Z$ such that $\ell(x) = 0$ and, for each $y \in T_d$, we have $\ell(z) = \ell(y) - 1$ for exactly one neighbor $z$ of $y$, while $\ell(z) = \ell(y) + 1$ for the other $d$ neighbors of $y$.  We can think of
$\ell(y)$ as being the generation number of $y$, with $x$ belonging to generation zero.

For $0 < \rho < 1$, define $$w_{\rho}(A) = \sum_{x \in A} \rho^{\ell(x)}.$$  For the basic contact process, we have
\begin{equation}
E[w_{\rho}(A_t')] \leq \exp[(\lambda d \rho + \lambda \rho^{-1} - 1) t] w_{\rho}(A_0').
\label{wineq}
\end{equation}
This can be seen from the differential equation in (4.4) on p. 80 of Liggett (1999).  Alternatively, equation (4.6) on p. 81 of Liggett (1999) gives (\ref{wineq}) with equality for the branching random walk, so the inequality follows from comparing the contact process with the branching random walk.

Since initially there is just one pathogen at $x$, we have $w_{\rho}(A_0') = 1$.  Also, since (\ref{wineq}) holds for all $\rho$, we can optimize the bound by choosing $\rho = 1/\sqrt{d}$.  If $\lambda < \frac{1}{2 \sqrt{d}}$, then $\gamma = \lambda d \rho + \lambda \rho^{-1} - 1 < 0$.  Since $w_{\rho}(A_t') \geq 1$
whenever $x \in A_t'$, it follows that $$P(x \in A_t') \leq E[w_{\rho}(A_t')] \leq \exp(\gamma t).$$  Since $\gamma < 0$, we have $$\int_0^{\infty} P(x \in A_t') \: dt < \infty.$$
Using Theorem 4 with $A = C = \{x\}$ to compare the basic contact process and the contact process with mutations, we get $P(x \in A_t) \leq P(x \in A_t')$ for all $t$.  Thus,
\begin{equation}
E \bigg[ \int_0^{\infty} 1_{\{x \in A_t\}} \: dt \bigg] = \int_0^{\infty} P(x \in A_t) \: dt < \infty.
\label{finexp}
\end{equation}
On the event that $x$ is infected infinitely many times, almost surely the Lebesgue measure of the set of infected times is infinite.  Therefore, (\ref{weakcond}) follows from (\ref{finexp}).
\end{proof} 

\begin{Rmk}
{\em Because $P(x \in A_t) \leq P(x \in A_t')$, we have
\begin{equation}
\lim_{t \rightarrow \infty} P(x \in A_t) = 0
\label{probto0}
\end{equation}
whenever this result holds for the basic contact process.  In particular, if $\lambda < \lambda_2$, where $\lambda_2$ is the upper critical value for the contact process, that is, the infimum over the set of $\lambda$ for which the contact process survives strongly, then (\ref{probto0}) holds for all $r \in [0, 1]$.  Pemantle (1992) showed that $\lambda_2 \geq 0.309$ when $d = 5$ and $\lambda_2 \geq 0.354$ when $d = 4$.  These bounds, combined with Theorem 1.1, show that when $d = 4$ or $d = 5$, there is
a range of $\lambda$ such that for all $r \in [0, 1]$, the contact process with mutations survives but (\ref{probto0}) holds.
For $d \geq 6$, Theorem 2 gives the stronger result of weak survival for a range of $\lambda$.}
\end{Rmk}

\begin{Rmk}
{\em Unlike the case of the contact process, many monotonicity statements for $A_t$ are
either false or hard to prove. To illustrate, we note that the following statement is false, even when
$S$ has only two elements:
\begin{center}
``For all $\lambda, t$ and $r$, $P^A(A_t\neq\emptyset)$ is increasing in $A$ for sets $A$
of a single type."
\end{center}
\noindent To see this, let $S=\{x,y\}$,
$f(A)=\lim_{\lambda\rightarrow\infty}P^A(A_t\neq\emptyset)$ for a fixed $t>0$,
and for simplicity take $r=1/2$. Then it is not hard to check that
$$f(\{x,y\})=e^{-t},\quad f(\{x\},\{y\})=(1+t)e^{-t},\quad f(\{x\})=\bigg(1+\frac t2\bigg)e^{-t}.$$
(Here $\{x,y\}$ is the configuration with two particles of the same type, $\{x\},\{y\}$
is the configuration with two particles of different types, and $\{x\}$ is the
configuration in which $x$ is occupied and $y$ is vacant.) Therefore,
$$f(\{x,y\})< f(\{x\})< f(\{x\},\{y\}).$$
In particular, we do not know whether the survival probability $P^{\{x\}}(A_t\neq\emptyset\ \forall t)$
is monotone in $r$ for fixed $\lambda$, even though Theorem 4 suggests that this
is probably the case.}
\label{monrmk}
\end{Rmk}

\section {Proofs of Theorems 1.1 and 1.2}

We start with the proof of Theorem 1.1.  The proof requires some extra care because of the monotonicity issues mentioned in Remark \ref{monrmk}.

To prove this result, we first establish that our particle system dominates a particle system that obeys the same rules except that only one birth can ever happen on any given site, and no births can occur on 
sites $y$ for which $\ell(y) < 0$, where $\ell$ is the function defined in the proof of Theorem 2.  Note that this second restriction amounts to running the process on the rooted $d$-ary tree with $x$ as the root, so that $x$ has $d$ neighbors and the other vertices have $d+1$ neighbors.  Therefore, it will suffice to give a condition for this modified particle system to survive.

To couple the two particle systems, we construct them from the same collection of Poisson processes.
For each nonzero integer $k$, let $(T^k_i)_{i=1}^{\infty}$ be the times of a Poisson process of rate $1$ on $[0, \infty)$.  For each pair of adjacent vertices $v$ and $w$ in $T_d$, let $(T^{v,w}_i)_{i=1}^{\infty}$ be the times of a Poisson process of rate $\lambda$ on $[0, \infty)$ and let $(\xi^{v,w}_i)_{i=1}^{\infty}$ be a sequence of i.i.d. Bernoulli random variables that are one with probability $r$ and zero with probability $1-r$.  All of these processes are assumed to be independent of one another.

To construct the two particle systems, we start with an individual of type one at $x$.  At the times $T^{v,w}_i$, if site $v$ is occupied and site $w$ is vacant, the individual at site $v$ gives birth at site $w$, except that for the second particle system, this birth is suppressed if there was a previous birth at 
site $w$ or if $\ell(w) < 0$.  If $\xi^{v,w}_i = 0$, then the new individual born at site $w$ has the same type as the individual at site $v$.  If $\xi^{v,w}_i = 1$, then this individual has a new type.  If this new type is the $k$th type that appears simultaneously in both particle systems, then we label it type $k$.  If it is the $k$th type to appear in one particle system but not the other, then we label the type $-k$.  Therefore, the types indexed by positive integers are coupled, while those indexed by negative integers are not.  At the times $T^k_i$, all individuals of type $k$ (if there are any) die.

Let $A_t$ be the set of occupied sites in the first particle system at time $t$, and let $B_t$ be the set of occupied sites in the second particle system.  For nonzero integers $k$, let $A_t^{k}$ be the set of sites at time $t$ inhabited by a type $k$ individual in the first process, and let $B_t^{k}$ be the set of sites at time $t$ inhabited by a type $k$ individual in the second process.  The following lemma is sufficient to imply that 
the first process survives for all $t$ whenever the second process does.

\begin{Lemma}
For all $t$, we have $B_t^k = \emptyset$ for all $k < 0$ and $B^k_t \subset A^k_t$ for all $k > 0$.
\label{couplelem}
\end{Lemma}

\begin{proof}
Fix any vertex $w \in T_d$.  Let $v_0, v_1, \dots, v_s$ be the vertices on the unique path from $x$ to $w$, with $v_0 = x$ and $v_s = w$.  Suppose $w \in B_t$ for some $t$.  This means that for $j = 1, \dots, s$, the virus must have spread from $v_{j-1}$ to $v_j$ by time $t$.  Therefore, $\ell(v_j) \geq 0$ for $j = 0, 1, \dots, s$.

If $T^1_1 < T^{v_0, v_1}_1$, then the individual at $v_0$ dies before the virus spreads to $v_1$.  In this case, the site $v_1$ can never become infected in the second process because, by the rules of the second process, $v_0$ can never be reinfected.  Alternatively, if $T^1_1 > T^{v_0, v_1}_1$, then the site $v_1$ becomes infected in both particle systems at time $\tau_1 = T^{v_0, v_1}_1$, and by definition the site is assigned the same (positive) type number in both particle systems.  Therefore, if $v_1$ is ever infected in the second particle system, then it is infected in both particle systems by particles of the same positive type at the same time $\tau_1$.

Proceeding by induction, suppose it is true that if $v_j$ is ever infected in the second particle system, then it is infected in both particle systems by a particle of the same type $m_j > 0$ at the same time $\tau_j$.  Then, if $\min\{T^{m_j}_i: T^{m_j}_i > \tau_j\} < \min\{T^{v_j, v_{j+1}}_i: T^{v_j, v_{j+1}}_i > \tau_j\}$, the individual at $v_j$ dies before spreading the infection to $v_{j+1}$, and the site $v_{j+1}$ 
can never become infected in the second particle system because $v_j$ can never be reinfected.  Otherwise, the site $v_{j+1}$ becomes infected in both particle systems at time $\tau_{j+1} = \min\{T^{v_j, v_{j+1}}_i: T^{v_j, v_{j+1}}_i > \tau_j\}$ and therefore gets the same type number in both particle systems.  Thus, if $v_{j+1}$ is ever infected in the second particle system, then it is infected in both particle systems at the same time by particles of the same positive type.

By induction, we conclude that for all $w \in T_d$, if $w$ is ever infected in the second particle system, then it is infected in both particle systems at the same time by particles having the same positive type.  This implies $B_t^k = \emptyset$ for all $k < 0$.  This also implies that the particle on site $w$ will die at the same time in both particle systems.  Since the site can never become reinfected in the second particle system, the result $B_t^k \subset A_t^k$ for all $k$ follows.
\end{proof}

\begin{Lemma}
Let $U$ be the tree such that $1$ is the root of the tree and $\ell$ is a child of $k$ if the first type $\ell$ individual in the second particle system has a type $k$ individual as its parent.  Then $U$ is a Galton-Watson tree.
\end{Lemma}

\begin{proof}
To prove this lemma, we will modify the Poisson process construction of the second particle system.  Rather than working with the Poisson processes $(T^{v,w}_i)_{i=1}^{\infty}$ and $(\xi^{v,w}_i)_{i=1}^{\infty}$, we will define, for each positive integer $k$, a Poisson process $(T^{k,v,w}_i)_{i=1}^{\infty}$ of rate $\lambda$ on $[0, \infty)$ and a sequence $(\xi^{k,v,w}_i)_{i=1}^{\infty}$ of i.i.d. random variables that are one with probability $r$ and zero with probability $1-r$.  All of the processes are again assumed to be independent.  We construct the particle system as before, except that when there is a type $k$ individual at the site $v$, the times at which this individual attempts to give birth on the site $w$ are the times $T^{k,v,w}_i$, and the random variables $\xi^{k,v,w}_i$ determine whether the new individual born has a new type or a different type.
Note also that is sufficient to consider only the vertices $v$ such that all vertices $w$ on the path from $x$ to $v$ satisfy $\ell(w) \geq 0$.  These vertices form a tree $T_d^*$, which we can think of as being rooted at $x$, in which $x$ has $d$ neighbors and the other vertices have $d+1$ neighbors.

For all positive integers $k$, let $X_k$ denote the number of children of the vertex labeled $k$ (if there is one) in the tree $U$.  Fix a vertex $u \in T_d$, and let $v_0, v_1, \dots, v_k$ be the vertices on the
unique path from $x$ to $u$.  Then, a type 1 individual gives birth to a new type on $u$ if and only if, for $j = 1, \dots, k-1$, the type 1 individual on site $v_{j-1}$ gives birth to another type 1 individual on site $v_j$, and then the type 1 individual on site $k-1$ gives birth to an individual of a new type on site $k$,
all before the type 1 individuals die.
These events depend only on $T_1^1$ and the processes $(T^{1,v,w}_i)_{i=1}^{\infty}$ and $(\xi^{1,v,w}_i)_{i=1}^{\infty}$, so there is a function $f$ such that $X_1 = f(T^1_1, (T^{1,v,w}_i)_{i=1}^{\infty}, (\xi^{1,v,w}_i)_{i=1}^{\infty})$.

Now, suppose there is a type $k$ individual born at some time.  Let $\tau_k$ be the time at which the first type $k$ individual is born, and let $v_k$ be the vertex at which the first type $k$ individual is born.  Note that $\tau_k$ and $v_k$ depend only on the evolution of the first $k-1$ types and therefore are 
functions of $T_1^j$, $(T^{j,v,w}_i)_{i=1}^{\infty}$, and $(\xi^{j,v,w}_i)_{i=1}^{\infty}$ for $j = 1, \dots, k-1$.  Also,
note that type $k$ individuals and their children can only be born in the subtree rooted at $v_k$ (consisting of the vertices $w$ such that the path from $w$ to $x$ includes $v_k$), which is empty until time $\tau_k$.  This means that the type $k$ individuals evolve in the subtree rooted at $v_k$ just as the
as the type $1$ individuals evolve in $T_d^*$.  More precisely, let $\phi$ be a graph automorphism mapping $T_d^*$ to the subtree of 
$T_d^*$ rooted at $v_k$ such that $\phi(x) = v_k$.
Also, let $J^k = \min\{i: T^k_i > \tau_k\}$ and $J^{k,v,w} = \min\{i: T^{k,v,w}_i > \tau_k\}$.  Then $$X_k = f \big( T^k_{J^k} - \tau_k, (T^{k,\phi(v),\phi(w)}_{i + J^{k,\phi(v),\phi(w)} - 1} - \tau_k)_{i=1}^{\infty}, (\xi^{k,\phi(v),\phi(w)}_{i + J^{k,\phi(v),\phi(w)} - 1})_{i=1}^{\infty} \big).$$
From this fact and the translation invariance of Poisson processes, it follows that the conditional distribution of $X_k$ given $X_1, \dots, X_{k-1}$ is the same as the unconditional distribution of $X_1$.  This is enough to establish that $U$ is a Galton-Watson tree.
\end{proof}

\begin{proof}[Proof of Theorem 1.1]
When $r = 1$, the result is a known fact about the basic contact process (see, for example, Theorem 4.1 on p. 79 of Liggett (1999)).  Assume now that $r < 1$.
By Lemma \ref{couplelem}, it suffices to prove that the second process survives with positive probability whenever $\lambda > 1/(d-1)$.  Note that the second process survives for all times if and only if the Galton-Watson tree $U$ is infinite.  Therefore, we need to show that this Galton-Watson process is supercritical when $\lambda > 1/(d-1)$

The mean of the offspring distribution can be computed explicitly.  We will take $k = 1$ and count the expected number of type 1 individuals born, not counting the one alive at time zero.
Note that for all $j \geq 1$, there are $d^j$ vertices $x \in T_d^*$ such that $\ell(x) = j$.  Also, if a type one vertex is born on a given site $x$, the probability that a type 1 vertex is born on a child of this site (that is, a neighboring vertex $y$ such that $\ell(y) = \ell(x) + 1$) is the probability that this birth, which happens at rate $\lambda (1-r)$, happens before either a mutant is born on this site, which happens at rate $\lambda r$, or the type dies, which happens at rate $1$.  Therefore, this probability is $\lambda(1-r)/(\lambda + 1)$.  It follows that the probability that a type $1$ individual is born on a given vertex $x \in T_d^*$ with $\ell(x) = j$ is $[\lambda(1-r)/(\lambda + 1)]^j$.  Therefore, the expected number of type $1$ individuals born is $$\sum_{j=1}^{\infty} d^j \bigg( \frac{\lambda(1-r)}{\lambda+1} \bigg)^j.$$  Since the expected number of new types born to type 1 individuals is $r/(1-r)$ times the expected number of type
one individuals born, the mean of the offspring distribution of the tree $U$ is
$$\frac{r}{1-r} \sum_{j=1}^{\infty} \bigg( \frac{d \lambda (1-r)}{\lambda + 1} \bigg)^j,$$
which is infinite if $r < (d-1)/d$ and $\lambda \geq 1/(d - 1 - dr)$.  Otherwise, we have $$\frac{r}{1-r} \sum_{j=1}^{\infty} \bigg( \frac{d \lambda (1-r)}{\lambda + 1} \bigg)^j = \frac{r}{1-r} \bigg( \frac{d \lambda (1-r)}{\lambda + 1 - d \lambda (1-r)} \bigg) = \frac{d \lambda r}{1 - (d-1)\lambda + d \lambda r},$$ which (when the mean of the offspring distribution is finite) is greater than one if and only if $\lambda > 1/(d-1)$.  This completes the proof of Theorem 1.1.
\end{proof}

We now turn to Theorem 1.2.  

\begin{proof}[Proof of Theorem 1.2]
When $r = 1$, the result is a known fact about the contact process.  By using a comparison between the contact process and the branching random walk, one can deduce this result from Theorem 4.8 in Section 1.4 of Liggett (1999).  When $r = 0$, all pathogens have the same type and therefore die out after an exponential time.  Assume now that $0 < r < 1$.

We first construct the original particle system from Poisson processes in a slightly different way.  For each positive integer $k$, let $(T_i^k)_{i=1}^{\infty}$ be the times of a rate one Poisson process on $[0, \infty)$.  For each pair of adjacent vertices $v$ and $w$ in $T_d$ and each positive integer $k$, let $(T^{k, v, w})_{i=1}^{\infty}$ be a Poisson process of rate $\lambda$ on $[0, \infty)$.  For each positive integer $k$ and each vertex $v$ in $T_d$, let $(\xi^{k,v})_{i=1}^{\infty}$ be a sequence of i.i.d. random variables that are one with probability $r$ and zero with probability $1-r$.
All type $k$ individuals die at the times $T_i^k$.  If, at time $T^{k,v,w}_i$, there is a type $k$ individual at site $v$ and the site $w$ is vacant, then the individual at site $v$ gives birth onto the site $w$ at this time.  To determine the type of this new individual, we examine the random variable $\xi^{k,w}_j$ if this is the $j$th time that a type $k$ individual has given birth onto $w$.  If $\xi^{k,w}_j = 0$, then the new individual also has type $k$.  Otherwise, it has a new type.  No types are assigned negative numbers, and the $\ell$th type to appear will be called type $\ell$.  Let $A^k_t$ be the set of sites having type $k$ at time $t$, and let $A_t$ be the set of occupied sites at time $t$.

To make our comparison argument work, we will define for each positive integer $k$ an additional process consisting only of type $k$ individuals.  We will define the new process from the same collection of Poisson processes as the original process.  If the first type $k$ individual in the original process is born at site $v_k$ at time $\tau_k$, then the new process will be empty until time $\tau_k$, at which time a type $k$ individual appears at the site $v_k$.  Thereafter, the new process evolves with the same rules as the original process, except that individuals of types other than $k$ are killed instantly after being born.  Let $D^k_t$ be the set of occupied sites in this process at time $t$.

Let $\zeta_k = \min\{T_i^k: T_i^k > \tau_k\}$, which is the time at which the type $k$ individuals die in both processes, so $A^k_t = D^k_t = \emptyset$ for $t \geq \zeta_k$.  For all $w \in T_d$, let $g(w) = \min\{j: \xi_j^{k,w} = 0\}$.
Note that for both processes, there is a type $k$ individual born at site $w$ if and only if there are at least $g(w)$ times before $\zeta_k$ when a type $k$ individual gives birth on site $w$.  Also, note that once a type $k$ individual is born on site $w$, there can be no further births on site $w$ until time $\zeta_k$, so there can never be more than $g(w)$ times at which a type $k$ individual gives birth on site $w$.  Let $h_{A,w}(t)$ be the number of times in the original process that a type $k$ individual gives birth on the site $w$ by time $t$, and let $h_{D,w}(t)$ be the number of times in the additional process that
an individual gives birth on the site $w$ by time $t$.

We claim that $h_{A,w}(t) \leq h_{D,w}(t)$ for all $t \in [\tau_k, \zeta_k]$ and all $w \in T_d$, and that $A_t^k \subset D_t^k$ for all $t \in [\tau_k, \zeta_k]$.  Since we clearly have $h_{A,w}(t) \leq h_{D,w}(t)$ and $A^t_k \subset D_t^k$ when $t = \tau_k$ and the type $k$ individuals die at the same time in both processes, we need only to show that these results do not break down at any of the times $T_i^{k,v,w}$ when type $k$ births may take place.  Suppose $t = T_i^{k,v,w}$, and that we have $h_{A,w}(t-) \leq h_{D,w}(t-)$ and $A_{t-}^k \subset D_{t-}^k$.  If $v \notin A_{t-}^k$, then there can be no birth in the original process at time $t$, and therefore $h_{A,w}(t) \leq h_{D,w}(t)$ and $A_{t}^k \subset D_{t}^k$.  
Suppose instead $v \in A_{t-}^k \subset D_{t-}^k$.  Then there are four possibilities.
\begin{itemize}
\item If $w \in A_{t-}$ and $w \in D^k_{t-}$, then there is no birth at time $t$ in either process, so
$h_{A,w}(t) \leq h_{D,w}(t)$ and $A_{t}^k \subset D_{t}^k$.

\item If $w \in A_{t-}$ and $w \notin D^k_{t-}$, then $h_{A,w}(t) = h_{A,w}(t-)$ and $h_{D,w}(t) = h_{D,w}(t-) + 1$.  Again it is clear that $h_{A,w}(t) \leq h_{D,w}(t)$ and $A_{t}^k \subset D_{t}^k$.

\item If $w \notin A_{t-}$ and $w \notin D^k_{t-}$, then there is a birth at time $t$ in both processes.  Therefore, $h_{A,w}(t) = h_{A,w}(t-) + 1$ and $h_{D,w}(t) = h_{D,w}(t-) + 1$, so we have $h_{A,w}(t) \leq h_{D,w}(t)$.  If the individual born in the original process has type $k$, then $h_{A,w}(t) = g(w)$, which implies $h_{D,w}(t) = g(w)$ and therefore a type $k$ individual is born in both processes.  Thus, we still have $A_{t}^k \subset D_{t}^k$.

\item If $w \notin A_{t-}$ and $w \in D^k_{t-}$, then it is clear that $A_{t}^k \subset D_{t}^k$.  Also, the individual at site $w$ at time $t-$ in the additional process must have type $k$, as only type $k$ individuals appear in this process.  This means that $h_{D,w}(t-) = g(w)$.  Since $h_{A,w}(t) \leq g(w)$ for all $t$, the inequality $h_{A,w}(t) \leq h_{D,w}(t)$ holds also.
\end{itemize}
These observations imply the claim.

Now, let $Y_k$ be the number of new types born to type $k$ individuals at the site $w$ in the original process, and let $Z_k$ be the number of new types born to type $k$ individuals at the site $w$ in the additional process.  Note that $Y_k = \min\{g(w) - 1, h_{A,w}(\zeta_k)\}$ and $Z_k = \min\{g(w) - 1, h_{D,w}(\zeta_k)\}$.  Therefore, it follows from the claim that $Y_k \leq Z_k$.  Also, because $T_d$ is a transitive graph, the conditional distribution of $Z_k$ given $Z_1, \dots, Z_{k-1}$ and given that a type $k$ individual appears at some time is the same as the distribution of $Z_1$.  

Using the original process, define the tree $V$ such that $1$ is the root of the tree and $\ell$ is a child of $k$ if the first type $\ell$ individual has a type $k$ individual as its parent.  Note that the vertex $k$ has $Y_k$ descendants in the tree $V$.  Since $Y_k \leq Z_k$ for all $k$, it follows that the distribution of the total progeny of $V$ is stochastically dominated by the distribution of the total progeny in a Galton-Watson tree whose offspring distribution is the same as the distribution of $Z_1$.  Therefore, if $E[Z_1] \leq 1$, then $V$ is finite almost surely, which implies that the process dies out almost surely.
 
The computation for $E[Z_1]$ differs from the computation in the proof of Theorem 1.1 in
two ways.  First, there are $(d+1) \cdot d^{j-1}$ vertices at a distance $j$ from $x$.  Second, because new types die immediately, if there is a type one present on a given site, a type one will be born on a child at rate $\lambda(1 - r)/(\lambda(1 - r) + 1)$ because this birth only must happen before the type dies.  Therefore, the expected number of new types born to type 1 individuals in this model is
$$\frac{(d+1)r}{d(1-r)} \sum_{j=1}^{\infty} \bigg( \frac{d \lambda (1-r)}{\lambda(1 - r) + 1} \bigg)^j.$$  This is infinite when $\lambda \geq 1/[(d-1)(1 - r)]$.  Otherwise, we have
$$\frac{(d+1)r}{d(1-r)} \sum_{j=1}^{\infty} \bigg( \frac{d \lambda (1-r)}{\lambda(1 - r) + 1} \bigg)^j =
\frac{(d+1)r}{d(1-r)} \bigg( \frac{d \lambda(1-r)}{\lambda(1-r) + 1 - d \lambda(1 - r)} \bigg)
= \frac{(d+1) r \lambda}{1 - (d-1)(1-r)\lambda},$$
which is less than or equal to one if and only if $\lambda \leq 1/(d - 1 + 2r)$.
\end{proof}

\section {Proof of Theorem 1.3}

In this section, we adapt the proof of Theorem 2.2(ii) of Pemantle (1992) to prove that the process
sometimes survives even for $\lambda < 1/(d-1)$. We begin with the following general result.

\begin{Lemma} Let $X_t$ be a pure jump Markov process on a countable set $S$ with
transition rates $q(x,y), x\neq y$, set $q(x)=\sum_{y: y\neq x}q(x,y),$ and
let $f$ be a nonnegative function on $S$. Suppose that
there exist positive $\epsilon$ and $M$ so that

(i) $f(y)-f(x)\geq -M$ whenever $q(x,y)>0$

\noindent and

(ii) $\sum_{y:|f(y)-f(x)|\leq M}q(x,y)[f(y)-f(x)]\geq\epsilon q(x)$ for all $x\in S$.

\noindent Then there is an $\alpha>0$ so that $E^x\exp[-\alpha f(X_t)]$ is a decreasing
function of $t$ for all $x\in S$. In particular, if $f(x)>0$, $\tau$ is the
hitting time of $\{y:f(y)=0\}$, and this set is absorbing, then $P^x(\tau<\infty)<1.$
\end{Lemma}

\begin {proof} Choose $\gamma>0$ so that $e^u\leq 1+u+u^2$ for $|u|\leq \gamma$, and write
for $0<\alpha\leq\gamma/M$

\begin{align}
\frac d{dt}E^xe^{-\alpha f(X_t)}\bigg|_{t=0} &= \sum_yq(x,y)\bigg[e^{-\alpha f(y)}-
e^{-\alpha f(x)}\bigg] \nonumber \\
&\leq e^{-\alpha f(x)}\sum_{y:|f(y)-f(x)|\leq M}q(x,y)
\bigg[e^{\alpha[f(x)-f(y)]}-1\bigg] \nonumber \\
&\leq e^{-\alpha f(x)}\sum_{y:y\neq x,|f(y)-f(x)|\leq M}
q(x,y)\bigg[\alpha[f(x)-f(y)]+\alpha^2M^2\bigg] \nonumber \\
&\leq e^{-\alpha f(x)}
\bigg[-\epsilon\alpha q(x)+\alpha^2M^2q(x)\bigg] \nonumber,
\end{align}
which is $\leq 0$ provided that $\alpha M^2\leq\epsilon$. Setting $\alpha=\min(\gamma/M,
\epsilon/M^2)$ and using the Markov property gives the monotonicity result. To check
the final statement, use the fact that if $P^x(\tau<\infty)=1$, then $f(X_t)$ is eventually
0. This would contradict the monotonicity statement.
\end{proof}

Returning to the mutation process, for a configuration $A$, let $N(A)$ be the number of
types in $A$, and $C(A)$ be the number of components of $A$. We will apply the lemma
with $f(A)=\alpha N(A)+\beta C(A)$ for an appropriate choice of $\alpha>0$ and $\beta>0$. Note that
$N(A)$ can only increase by 1 or decrease by 1 at each transition. On the other hand,
$C(A)$ can decrease by at most $d$ (if a particle is born at a site with all neighbors
occupied), but can increase by arbitrarily large amounts (if a type that is connected to
a large number of other types dies). So $f(A_t)$ has bounded downward jumps but unbounded
upward jumps, unlike the situation Pemantle considered. Nevertheless, by choosing
$M$ sufficiently large, which we now do, we satisfy (i) of the lemma.

Turning to assumption (ii), we consider first the contributions coming from changes in $N(A)$.
The number of types decreases by one at rate $N(A)$, since each type dies at rate 1.
The number of types increases by one at rate
$$\lambda r\#((x,y):x\in A, y\notin A)=\lambda r[(d-1)|A|+2C(A)].$$
This is exactly the same as Pemantle's computation, except that there is an extra factor
of $r$ coming from the fact that each new particle is of a new type with probability $r$.
Combining these observations leads to
$$\sum_Bq(A,B)[N(B)-N(A)]=-N(A)+\lambda r[(d-1)|A|+2C(A)].$$

The rate at which $C(A)$ decreases as the result of particle
births is at most $(d+1)\lambda C(A)$, as shown by Pemantle. The major difference
between Pemantle's situation and ours comes next. Consider one component $C$ of $A$ with $k\geq 2$ types,
in which type $i$ is connected to $j_i$ other types. As Pemantle showed,
\begin{equation}
\sum_{i=1}^k(j_i-1)=k-2.
\label{oldeq2}
\end{equation}
If the $i$th type dies, $j_i-1$ new components are formed. Thus we need to find a lower bound
for
$$\sum_{i:j_i\leq M+1}(j_i-1).$$
In doing so, we will use (\ref{oldeq2}), together with the simple inequalities
\begin{equation}
j_i\leq (d-1)|A^i|+2,
\label{oldeq3}
\end{equation}
\begin{equation}
|C|=\sum_{i=1}^k|A^i|\geq\sum_{i:j_i> M+1}|A^i|+\sum_{i:j_i\leq M+1}1
\label{oldeq4}
\end{equation}
and (using (\ref{oldeq2}) and the fact that $j_i\geq 1$)
\begin{equation}
k-2\geq\sum_{i:j_i> M+1}(j_i-1)\geq (M+1)\sum_{i:j_i> M+1}1.
\label{oldeq5}
\end{equation}
Putting these together gives:

\begin{align}
\sum_{i:j_i\leq M+1}(j_i-1) &= (k-2)-\sum_{i:j_i> M+1}(j_i-1) \nonumber \\
&\geq(k-2)-\sum_{i:j_i>
M+1}\big[(d-1)|A^i|+1\big] \nonumber \\
&\geq (k-2)-\sum_{i:j_i> M+1}1-(d-1)\bigg[|C|-\sum_{i:j_i\leq
M+1}1\bigg] \nonumber \\
&=(k-2)-(d-1)|C|+(d-1)k-d\sum_{i:j_i>
M+1}1\nonumber \\
&\geq kd-2-(d-1)|C|-d\frac{k-2}{M+1} \nonumber \\
&\geq kd\frac{M}{M+1}-(d-1)|C|-2. \nonumber
\end{align}

\noindent We have used (\ref{oldeq3}) in the first inequality, (\ref{oldeq4}) in the second, and (\ref{oldeq5}) in the third. Note that this
inequality is trivially true in the case $k=1$ in which the entire component is of one type.
Summing the above inequality over all components gives
$$\sum_{i=1}^{N(A)}(j_i-1)1_{\{j_i\leq M+1\}}\geq N(A)d\frac {M}{M+1}-(d-1)|A|-2C(A).$$
It follows that
$$\sum_{B:|C(B)-C(A)|\leq M} q(A, B)[C(B)-C(A)]\geq N(A)d\frac M{M+1}-(d-1)|A|-2C(A)-(d+1)\lambda C(A).$$

Now take $f(A)=\alpha N(A)+\beta C(A)$.  Note that when $M$ is sufficiently large, we can only have $|f(B) - f(A)| > M$ with $q(A, B) > 0$ when the death of a single type leads to the creation of many components.  That is, $|f(B) - f(A)| > M$ if and only if $N(B) - N(A) = -1$ and $C(B) - C(A) > (M + \alpha)/\beta$.  It follows that for sufficiently large $M$,
\begin{align}
&\hspace{-.1in}\sum_{B:|f(B)-f(A)|\leq M}q(A,B)[f(B)-f(A)] \nonumber \\
&\geq
\sum_{B} q(A,B)[N(B) - N(A)] + \sum_{B:|C(B)-C(A)|\leq M/\beta} q(A,B)[C(B)-C(A)] \nonumber \\
&\geq \alpha\bigg[-N(A)+\lambda r(d-1)|A|+2\lambda
rC(A)\bigg] \nonumber \\
&\hspace{.2in}+\beta\bigg[N(A) d\frac{M}{M+\beta}-(d-1)|A|-2C(A)
-(d+1)\lambda C(A)\bigg] \nonumber \\
&=N(A)\bigg[-\alpha+\beta d\frac{M}{M+\beta}\bigg]+|A|\bigg[\alpha\lambda r(d-1)-\beta
(d-1)\bigg]+C(A)\bigg[2\lambda
r\alpha-2\beta-\beta(d+1)\lambda\bigg]. \nonumber
\end{align}

Since the total jump rate $q(A)$ is bounded above by a multiple of $|A|$,
and $N(A)\leq |A|$, assumption (ii) in
the lemma will be satisfied if, in the expression above, the coefficient of $|A|$ is strictly
positive, the sum of the coefficients of $N(A)$ and $|A|$ is strictly positive, and the
coefficient of $C(A)$ is nonnegative. This will be true for large $M$ if
$$\alpha\lambda r>\beta,\quad \beta>\alpha[1-\lambda r(d-1)],\text{ and }2\lambda r\alpha\geq
\beta[2+\lambda(d+1)].$$
There exist positive $\alpha$ and $\beta$ satisfying these inequalities if $\lambda r(d-1)\geq1$,
or if $\lambda r(d-1)<1$ and
$$\frac{2+\lambda(d+1)}{2\lambda r}<\frac 1{1-\lambda r(d-1)},$$
i.e., if $\lambda$ is greater than the positive root of
$$\lambda^2r(d^2-1)+\lambda(2rd-d-1)-2=0.$$
Applying the lemma, we now have the following result.

\begin{Theo} The contact process with mutations survives if
\begin{equation}
\lambda>\frac{d+1-2rd+\sqrt{(d+1)^2+4r(d+1)(d-2)+4d^2r^2}}{2r(d^2-1)}.
\label{lambdabound}
\end{equation}
\end{Theo}

Combining this result with Theorem 1.2, we see that for
$$\frac{1-d+\sqrt{(d-1)(7+9d)}}{2(d^2-1)}<\lambda <\frac 1{d-1},$$
the process dies out for small $r$ and survives for $r$ close to 1.
Note also that while Theorem 11 gives a better lower bound than Theorem 1.1 when $r$ is close to one, the lower bound in Theorem 1.1 is better for small $r$, as the expression on the right-hand side of (\ref{lambdabound}) approaches infinity as $r \rightarrow 0$.

\bigskip
\bigskip
\bigskip
\begin{center}
{\bf {\Large References}}
\end{center}

\mn N. Berger, C. Borgs, J. Chayes, and A. Saberi (2005).  On the spread of viruses on the internet.  In {\it Proceedings of the 16th ACM-SIAM Symposium on Discrete Algorithms}, 301-310.

\mn M. Draief, A. Ganesh and L. Massouli\'e (2006). Thresholds for virus spread on networks. Preprint, arXiv:math.PR/0606514.

\mn A. Ganesh, L. Massouli\'e, and D. Towsley (2005).  The effect of network topology on the spread of epidemics.  In {\it Proceedings of IEEE INFOCOM 2005}, 1455--1466.

\mn T. E. Harris (1974). Contact interactions on a lattice. {\it Ann. Probab.}, {\bf  2}, 969--988.

\mn T. M. Liggett (1985). {\it Interacting Particle Systems.} Springer, New York.

\mn T. M. Liggett (1999). {\it Stochastic Interacting Systems: Contact, Voter and Exclusion Processes}, Springer, Berlin.

\mn R. Pastor-Satorras and A. Vespignani (2001).  Epidemic spreading in scale-free networks.  {\it Phys. Rev. E}, {\bf 86}, 3200-3203.

\mn R. Pemantle (1992) The contact process on trees. {\it Ann. Probab}, {\bf 20}, 2089--2116.

\mn R. B. Schinazi and J. Schweinsberg (2006). Spatial and non-spatial stochastic models for immune response. To appear in {\it Markov Process. Related Fields}.
\end{document}